\documentclass[10pt,reqno]{article}

\newif\iffigures\figurestrue
\figuresfalse

\usepackage[colorlinks=false,pdfborderstyle={/W 1}]{hyperref}
\newif\ifhyper\IfFileExists{hyperref.sty}{\hypertrue}{\hyperfalse}

\ifhyper\usepackage{hyperref}
\def\hitem#1#2{\item[\hypertarget{#1}{#2}]\expandafter\gdef\csname LBL#1ITM\endcsname{#2}}
\def\iref#1{\hyperlink{#1}{\csname LBL#1ITM\endcsname}}
\else
\def\hitem#1#2{\item[{#2}]\expandafter\gdef\csname LBL#1ITM\endcsname{#2}}
\def\iref#1{{\csname LBL#1ITM\endcsname}}
\fi

\usepackage{amssymb,amsmath,amsthm,amsfonts,mathrsfs,ulem, cancel}
\usepackage[utf8]{inputenc}

\def\F{\mathscr{F}}

\def\text#1{\hbox{#1}}

\def\N{{\mathbb N}}
\def\T{{\mathbb T}}

\def\P{\mathbf P}
\def\E{\mathbf E}

\def\1{\mathbf{1}}

\newtheorem{theorem}{Theorem}[section]
\newtheorem{definition}{Definition}[section]
\newtheorem{lemma}[theorem]{Lemma}

\newtheorem{remark}{{\bf Remark}}
\newtheorem*{note}{\underline{Note}}
\newtheorem{example}{{Example}}[section]

\long\def\note#1/{\ifdraft{\marginpar{{$\Longleftarrow$}} \bf [#1] }\fi}

\newcounter{my}
\stepcounter{my}

\def\1{1\!\! 1}

\def\FIID{\mathsf{FIID}}
\def\clust{\mathsf{clust}}
\def\Ber{\mathsf{Ber}}

\def\un{\mathsf{un}}

\def\Vor{\mathsf{Vor}}

\theoremstyle{definition}
\numberwithin{equation}{section}
\numberwithin{figure}{section}

\title{A factor of i.i.d.\ with uniform marginals
and infinite clusters spanned by equal labels}

\author{Péter Mester\footnote{Alfréd Rényi Institute of Mathematics} 
}
\date{ \today}

\begin{document}

\maketitle

\begin{abstract}
We give an example of an $\FIID$ vertex-labeling of $\T_3$ whose marginals are uniform on $[0,1]$, and if we delete the edges between those vertices whose labels are different, then some of the remaining clusters are infinite.
\end{abstract}

\section{Introduction}

In this paper we answer a question asked by Gaboriau in his ICM survey in a slightly different language {(Question 5.6 in \cite{Gab})}. The question concerns factor of i.i.d.\ processes ($\FIID$s from now on), which are random elements of ${S^{\Gamma}}$ where $\Gamma$ is a countable group and $S$ is some {\it label set} (for us it will be $[0,1]$ or $\{0,1\}^n$), and arise by applying an equivariant measurable map $f:[0,1]^{\Gamma}\rightarrow S^{\Gamma}$ to an i.i.d.\ family $\{\omega_{\gamma }\}_{\gamma\in \Gamma}$ where $\omega_{\gamma}$ is uniform on $[0,1]$. If $\Theta=f(\omega)$ is an $\FIID$ and $\gamma\in \Gamma$, then $\Theta(\gamma)$ is a random variable whose distribution (which does not depend on $\gamma$ by the equivariance of $f$) will be called the {\it marginal} of $\Theta$.
For us the group will be relevant as the vertex set of its Cayley graphs and we will give further intuition on $\FIID$s in that case.

If ${\cal L}\in S^V$ is a labeling of the vertex set $V$ of a graph, then let ${\clust}({\cal L})$ be the subgraph obtained by deleting any edge whose endpoints got different labels by ${\cal L}$. The connected components will be called {\it clusters}.

Gaboriau asked if assuming that the marginals of an $\FIID$ on a Cayley graph are uniform on $[0,1]$ implies that the corresponding clusters are finite? We will show by an example on the $3$-regular tree $\T_3$ that the answer is {\it no}:

\begin{theorem}[]\label{t.gab}

There is an $\FIID$ labeling $\Theta$ of $\T_3$, whose marginals are uniform on $[0,1]$ for which $\clust(\Theta)$ contains infinite clusters.

\end{theorem}

Actually $\T_3$ can be replaced by any nonamenable Cayley graph as \cite{Russ_Damien} has shown that they always contain an $\FIID$ spanning forest whose components have furcation vertices (defined in ~Section \ref{s.furcations}) and that is the only thing our construction needs.

We will actually provide two slightly different constructions, where the first one can be said to be more natural and the second one adds a little twist to it so that it will have the additional property that {\it every} cluster will be infinite. 

It is natural to think of $\FIID$s as the class of random labelings of the vertex set which can be obtained by applying a local relabeling algorithm ($f$ in the definition) which does not use any further randomness beyond its input ($\omega$ in the definition). 
The reason we call this algorithm local is that the measurablility of $f$ implies that knowing the restriction of $\omega$ to a large enough ball $B(v,r)$ allows us to approximate $f(\omega)(v)$, as the conditional expectation $\E \big( f(\omega)\, \big|\,\omega_{B(v,r)}\big)$ converges 
to $f(\omega)$ as $r\rightarrow\infty$. See \cite{RussFIID} for a survey with more references on $\FIID$s on trees.

In general it seems hard to decide whether or not an $\FIID$ process with a given property exists. The locality mentioned above suggests that the answer to Gaboriau's question should be affirmative since the conditions imply that the vertices build a label whose specific value have probability zero, through a process governed by local data so it seems surprising that vertices constituting whole infinite clusters could end up find the same label. There is in particular a strong correlation decay found by \cite{BSzV} which can be an obstruction for some processes to be $\FIID$ which might be seen as a quantitative version of this locality. While we will show that it is possible, there is an intuition that the condition should imply that at least the clusters are {\it small} in some sense. And indeed, a result by \cite{Chifan_Ioana} implies that the clusters under the condition that the marginals are uniform on $[0,1]$ must be {\it hyperfinite}. In our context hyperfiniteness of the forest ${\clust}(\Theta)$ means that there is a sequence of random invariant forests $\{F_i\}_{i\in \N}$ such that every component of $F_i$ is finite almost surely and (using the notation $F(v)$ for the component of $F$ containing the vertex $v$ in a forest $F$) $F_i(v)\subset F_j(v)$ for $i\leq j$ and ${\clust}(\Theta)(v)=\bigcup_{i\in \N}F_i(v)$ for each vertex $v$. While it is itself true for any countable set that it is an increasing union of finite sets, the extra requirement of achieving it with random invariant proecesses makes hyperfiniteness a strong property. It implies that the clusters of ${\clust}(\Theta)$ must be finite or must have $1$ or $2$ ends (see \cite{BLPS}).

The following problem will not be pursued in this paper, but the example serves as a quick illustration of this locality phenomenon of $\FIID$s contrasted with invariant processes in general. 

\begin{example}\label{x.Independent_Vertices} Assume we want an invariant random vertex set $S$ of $\T_d$ which is  {\it independent} in the graph theoretic sense (that is, $S$ does not contain neighbors). We want the marginal probability $p_S$ of being in $S$ to be high. The optimal $p_S={1 \over 2}$ is achieved by the following random invariant process: Let $v_1, v_2 \in V(\T_d)$ be equivalent if their distance is even, and pick one of the classes with probability ${1 \over 2}$ to be $S$.
If we want $S$ to be a $\FIID$ (i.e., its indicator $1_S$ to be  $\FIID$), then a possible solution is the following code: at a vertex $v$ let $f(v):=1$ iff for each neighbor $w$ of $v$, $\omega(v)>\omega(w)$ (where $\omega$ was the original source) and let $f(v)=0$ otherwise, then $S=\{f(v)=1:v\in V(\T_d)\}$ is an independent set. In this case $p_S={1\over{d+1}}$.

 While this simple construction can be improved, its basic features are known to hold even for a near-optimal independent set arising as an $\FIID$. Namely, the marginal of any $\FIID$ independent vertex set of $V(\T_d)$ is bounded away from ${1 \over 2}$ for any $d\geq 3$ and it goes to zero as $d\rightarrow \infty$.
This follows easily from \cite{Bollobas}, however, its focus is on finite graphs. For the connection with $\FIID$s and references to further research in this direction, see \cite{RV17}.
\end{example}

An example similar to the above will be used in our construction.

\begin{example}\label{x.Sparse_Vertices} Let $r$ be a positive integer and let $S$ be the $\FIID$ vertex set defined by the following code (which will be the indicator of $S$): let $1_S(v)=1$ iff the label $\omega(v)$ of $v$ is maximal of all the labels within the ball of radius $r$ around $v$ (otherwise $1_S(v)=0$). Then $S$ has the property that any two vertices of it have distance at least $r$. 
\end{example}

We will need this construction not only directly on $\T_3$  but also on some locally finite forests associated to it. Then ``distance'' will refer to the distance within the forest (thus infinite between vertices of different components). These forests will be random and $S$ can be sampled independently of it. Notice that in this way $S$ guaranteed to intersect all infinite component of these forests (there will only be countably many components).

We close this section with an important elementary observation.

By a ${\un}[0,1]$ random
variable we mean one which is uniform on $[0,1]$.  By $x\buildrel{d}\over\sim y$ we mean that the random variable $x$ and $y$ has the same distribution, but with some abuse of notation we will also denote by $x\buildrel{d}\over\sim \nu$ if $x$ has distribution $\nu$.
From a single $x\buildrel{d}\over\sim {\un} [0,1]$ we can obtain an i.i.d.\ family of infinitely many $x_j\buildrel{d}\over\sim {\un} [0,1]$ by reorganizing the bits of $x$. Using this when we describe the code $f$, we can assume that it can always reach out for an additional ${\un}[0,1]$ random variable independent of any other step of the algorithm. But, importantly, every single random variable the algorithm uses is local in the sense that it must belong to some vertex.

\section{Voronoi partitions and other forests}\label{s.Voronoi}

So far we have defined $\FIID$s only as process that are labelings of the vertices, we extend the notion to processes which are $\{0,1\}$-labelings of the edges, so that then they can be used to encode subgraphs (which are subforests in our case). The notion of equivariance goes through as an action of a group $\Gamma$ on $V$ extends naturally to $V\times V$.
For example, we have already defined the forest $\clust(L)$ corresponding to a labeling $L$.

We now define Voronoi partitions. 
If $S$ is a vertex set, we want to partition all the other vertices into classes according to the closest element of $S$.  We have to deal with the potential ambiguity if a vertex $v$ is at an equal distance from several elements of $S$, moreover, we want to make the partition classes connected. 

\begin{definition}\label{d.Voronoi} Let $F$ be any locally finite forest and let $S\subset V(F)$ and a collection of distinct real numbers $\{\alpha(v)\}_{v\in S}$ be given. If $v\in V(F)$, and the $F$-component which contains $v$ also contains some element from $S$, then let $S_v\subset S$ be the set of those elements
of $S$ which are closest to $v$, i.e., $S_v:=\{s: {
d}_F(v,s)={d}_F(v,S)\}$.  Let
$\phi_{(S,\alpha^{})}(v):=s_0$ be that element of $S_v$ for which
$\alpha^{}(s_0)$ is minimal (by the local finiteness of $F$, $S_v$ is finite). Let two vertices $v_1,v_2$ be
equivalent if $\phi_{(S,\alpha^{})}(v_1)=\phi_{(S,\alpha^{})}(v_2)$. 
If the $F$-component of $v$ does not contain any element from $S$, then let the equivalence class of $v$ be the singleton $\{v\}$.
Let ${\Vor}(S,\alpha^{})$ be the partition corresponding to
this equivalence. 
 \end{definition}

The role of the $\alpha$ is to handle the ambiguity if $|S_v|>1$, in this way the partition classes are indeed  connected (note that it is not true for just any convention that breaks the tie). 

When we use Voronoi partitions, the forest will be in the form of $\clust(L)$ or something closely related, the $S$ will be an $\FIID$ set and the $\alpha$ will be extracted from the source. We will suppress $\alpha$ in the notation and just denote the partition by ${\Vor}(S)$. We always assume that the hidden $\alpha$ is independent of any other steps of the construction. We will refer to partition classes as {\it cells} and we will use this terminology in general where we have a forest where every component is finite almost surely (we will see that Voronoi partitions have this property).
When we consider a Voronoi partition as a forest, we delete edges between vertices of different cells and forget the distinguished vertex; in this way they are $\FIID$ forests.

If we want to produce an $\FIID$ labeling with a $\un[0,1]$ marginal whose clusters are finite but arbitrarily large, that is easy. We can even sample an arbitrary random forest whose components are almost surely finite, and label the vertices independently afterwards:

\begin{lemma}\label{l.Finite_Clusters_Easy}
Let $\Pi$ be an ${\FIID}$ forest whose components are almost surely finite. There is an $\FIID$ labeling $\theta$ with a $\un [0,1]$ marginal which is constant over each component of $\Pi$ (in fact, almost surely $\clust(\theta)=\Pi$), and the $\theta$-labels of different components of $\Pi$ form an independent family.

\end{lemma}

\proof Let  $({\alpha}(v),\beta(v))_{v\in V(T)}$ be a collection of two independent ${\un}[0,1]$ label over each vertex. For a vertex $v$ let $\Pi(v)$ be the component of $\Pi$ containing $v$. Since almost surely $|\Pi(v)|<\infty$ and the $\beta(v)$ labels are all distinct, there will be a unique $v_0\in \Pi(v)$ for which $\beta(v_0)$ is minimal within $\Pi(v)$ (i.e.,  $\beta(v_0)={\tt
min}\{\beta(w);w\in \Pi(v)\}$.) Then let $v$ ``copy" the $\alpha$ label from $v_0$, meaning that $\theta(v):={\alpha}(v_0)$.
  \qed
  
  \medskip
  
  The finiteness of the components above was crucial; when we construct the infinite clusters with uniform labels, then the labels and the clusters will be built together step-by-step and not by selecting the infinite clusters first and labeling them afterwards.

We will call a random forest whose components are almost surely finite a {\it cell-partition} and the components will be called {\it cells} or $\Pi$-cells where the forest is denoted by $\Pi$.
If $\Pi$ is an $\FIID$ cell-partition, then let $\Ber(\Pi)$ be the $\FIID$ $\{0,1\}$-labeling $\lambda$ with the properties that: its marginals are fair bits ($\P(\lambda(v)=1)=\P(\lambda(v)=0)={1\over 2}$), the labels are constant over a $\Pi$-cell ($\Pi(v_1)=\Pi(v_2)$ implies $\lambda(v_1)=\lambda(v_2)$) and the labels over different cells are independent.  By the notation $\lambda\sim \Ber(\Pi)$ we will mean that first $\Pi$ is sampled, and then (given $\Pi$) $\lambda$. The notation refers to Bernoulli site percolations and we will use $\Ber(\Pi)$ to imitate a Bernoulli percolation on a graph whose vertices are the cells of $\Pi$. The fact that a labeling with this distribution can be realized as an $\FIID$ labeling is a consequence of the previous lemma: the fair bits needed for $\Ber(\Pi)$ can be obtained from the $\un[0,1]$ $\theta$-label guaranteed by the lemma, for example by defining the bit to be $1$ if $\theta>{1\over 2}$ and $0$ if $\theta\leq {1\over 2}$.

One may notice that for a Voronoi type partition ${\Vor}(S)$ we do not need to know the finiteness of the cells to label them as claimed in the lemma since each cell already comes with a single distinguished vertex (the one from $S$) and the whole cell can copy labels from this distinguished vertex just as in the proof.

So, if a Voronoi partition had nonzero chance of producing infinite clusters, then that already would witness the truth of our ~Theorem \ref{t.gab}.
However, a simple application of the Mass Transport Principle (see ~Chapter 8 in \cite{LPbook}) shows that every Voronoi cell must be finite.

\begin{lemma}\label{l.MassTransport} If an invariant process ${\cal R}$ on a Cayley graph produces connected components with a single distinguished vertex from each component, then each of these connected components must be finite. 
\end{lemma}

Before the proof, note that this immediately implies that this is also true if the word ``single'' is replaced by finitely many, as from the finitely many vertices we can select a uniform one and this still will be an invariant process if the original one was.

\proof Define the following function (which is invariant under graph automorphisms): $F(x,y,\omega):=1$ if $x$ is the distinguished vertex of the partition
class containing $y$ in the random configuration $\omega$. We will call $F(x,y, \omega)$ {\it the mass sent by
$x$ to $y$} or {\it the  mass received by $y$ from $x$}.
The Mass-Transport Principle says that if ${\cal R}$ is invariant,
then for the identity $o \in V$ the expected overall mass $o$ receives
is the same as the expected overall mass it sends out.
If there was a counterexample to the statement of the lemma, then the expected mass the origin would receive would be no more than
one (this is true even pointwise). However, the expected mass
it would send out is infinite (it even would send out infinite mass
with positive probability). \qed
\medskip

Recall that if $\Pi$ is a forest, then for a vertex $v$ we denoted by $\Pi(v)$ the component of $v$. If two forests $P,F$ are related in a way that $P(v)\subset F(v)$ for all $v$, then we denote this relationship by $P\prec F$ or $F\succ P$. To such a pair we associate a new forest:

\begin{definition}\label{d.LargeScale} If $F,P$ are forests on the same vertex set  and $P\prec F$, then we associate to this pair a new forest $F/P$ called the {\bf large scale forest} (or when $F$ is a tree, the {\it large scale tree}). The vertices of $F/P$ are the components of $P$ and two $P$-components $t_1,t_2$ are connected in $F/P$ if their distance is $1$ in $F$. For a vertex $v$ let $F/P(v)$ be the subtree of $F/P$ which contains $P(v)$.
\end{definition}

When we use this large scale forest construction $P$-components will be finite (so $P$ is a cell-partition). If there is a further cell-partition $\Pi$ on $F/P$, then there is a natural corresponding cell-partition ${\tt glue}_{\Pi}(P)$ on $F$ so that $P \prec {\tt glue}_{\Pi}(P) \prec F$. We just glue together the cells of $P$ according to $\Pi$, meaning that if ${\cal C}$ is a $\Pi$-cell consisting of the $P$-cells $C_1,\dots,C_l$, then $\bigcup\{C_i:C_i\in {\cal C}\}$ will be a ${\tt glue}_{\Pi}(P)$-cell and, as these cells already partition all the vertices of $F$, defines ${\tt glue}_{\Pi}(P)$.

Note also that, in the case $F,P$ are ${\FIID}$ subforests of $\T_3$ and $P\prec F$ and $P$ is a cell-partition, then by ~Lemma \ref{l.Finite_Clusters_Easy} we can assume that the vertices of $F/P$ are equipped with a family of i.i.d.\ random variables $x(v)_{v\in V(F/P)}\buildrel{d}\over\sim {\un} [0,1]$ which we can use to build Voronoi type partitions on $F/P$ as an ${\FIID}$-forest on the original $\T_3$.

\section{High Level Overview}\label{s.High_Level}

To highlight the ideas of the construction, we first show the modest claim that for any positive integer $n$, there is an $\FIID$ labeling $\theta_n$ on $\T_3$ whose marginal is uniform on the label set $\{0,1\}^n$ and $\clust(\theta_n)$ contains infinite clusters. We will use the basic theory of Bernoulli percolation on trees, see \cite{LPbook} or \cite{PGG}.

A Bernoulli-${p}$ site percolation (${\Ber}({p})$-labeling from now on) is the labeling where each vertex is labeled with $1$ with probability $p$ and with $0$ with probability $1-p$ independently of the others. On $\T_d$ it has infinite clusters whose label is constant $1$ exactly if $p(d-1)>1$.  If $L_0,\dots, L_{n-1}$ are independent $\Ber({1\over 2})$-labelings, and we concatenate them to get the $\{0,1\}^n$-labeling ${\cal L}_n:=(L_0,\dots,L_{n-1})$, then for any $s\in \{0,1\}^n$ the distribution of vertices whose ${\cal L}_n$-label is $s$ will be the same as the distribution of vertices whose label is $1$ in a single $\Ber({1\over {2^n}})$-labeling. In particular, there will be infinite clusters in $\clust({\cal L}_n)$ on $\T_d$ if $d>2^n+1$. 

This is not yet the $\theta_n$ we promised, as we want to label $\T_3$ instead of $\T_d$ where $d$ depends on $n$. But we can imitate a tree whose minimal degree is at least $d$ within $\T_3$ using Voronoi partitions and the large scale tree construction. As in ~Example \ref{x.Sparse_Vertices}, let $S$ be an $\FIID$ vertex set in $\T_3$ for which any two $v_1,v_2\in S$ has distance at least $2r+1$, then let $\Pi_0:=\Vor(S)$. Each $\Pi_0$-cell contains the ball $B_{\T_3}(v,r)$ around each $v\in S$. This implies that the large scale tree $\T_3/\Pi_0$ has minimal degree at least $|B_{\T_3}(v,r)|+2$.

 For large enough $r$ the random tree $\T_3/\Pi_0$ has minimal degree at least $2^n+2$, so using the $\FIID$ labeling $\Ber(\Pi_0)$ guaranteed by Lemma \ref{l.Finite_Clusters_Easy} and concatenating $n$ independent versions $\lambda_i\sim \Ber(\Pi_0)$ to form $\theta_n$, we get an $\FIID$ $\{0,1\}^n$-labeling which has infinite clusters and marginals uniform on $\{0,1\}^n$.

This proves the claim, but how do we get a labeling which has uniform marginals on $[0,1]$? If we keep adding extra independent bits to the already constructed $\theta_n$ by further $\lambda_j\sim\Ber(\Pi_0)$-labelings and take the sequence of bits as the binary representation of a real from $[0,1]$, then we would get a labeling with $\un[0,1]$-marginals.  Of course in this labeling there would not be any infinite clusters.

However, if instead of the ``static'' sequence $\Pi_0,\dots,\Pi_0,\dots,$ with the i.i.d.\ labels $\lambda_i\sim \Ber(\Pi_0)$, we use a dynamically changing sequence of cell-partitions $\Pi_0\prec\dots\prec\Pi_n\prec\dots$, and the corresponding sequence of labelings $\Lambda_0\sim \Ber(\Pi_0),\dots,\Lambda_n\sim\Ber(\Pi_n),\dots$, then we will be able to use infinitely many bits and thus getting $\un[0,1]$-marginals, while also having infinite clusters.  The essence of how this sequence is constructed and what issue needs to be taken care of is already visible in the step from $\Pi_0$ to $\Pi_1$. 

There are ``target degrees'' $D_0$ and $D_1$ which for now are just large integers (a target degree was implicit before where $\T_3/\Pi_0$ had minimal degree at least $2^n+2$).
$\Pi_0$ and $\Lambda_0\sim\Ber(\Pi_0))$ are defined as before, where now we want $\T_3/\Pi_0$ to have minimal degree at least $D_0$. We get the random forest ${\cal F}_0:=\clust(\Lambda_0)$, and of course $\Pi_0\prec{\cal F}_0$. We want to build $\Pi_1$ in such a way that $\Pi_0\prec\Pi_1\prec{\cal F}_0$ and ``whenever possible'' the components of the large scale forest ${\cal F}_0/\Pi_1$ should have minimal degree at least $D_1$. So the goal in this second step is similar to the one in the first step, when we wanted $\T_3/\Pi_0$ to have minimal degree at least $D_0$.

A key difference is that in the first step we worked in the the known tree $\T_3$, while now we have to deal with the random forest ${\cal F}_0$. We can see immediately that the target degree goal cannot be reached for all components of ${\cal F}_0/\Pi_1$, as ${\cal F}_0$ contains finite clusters. One can also build infinite trees in an adversarial way (see ~Example \ref{x.PathologicalFurcation}) which are obstacles to this goal. However, ${\cal F}_0$ is defined by an invariant random process which avoids those sorts of examples (by the same application of the Mass Transport Principle we saw before). But a random invariant process may produce a bi-infinite path, which would also be an obstacle to our goal. The kind of random components we will need are the ones which contain furcation vertices (defined in the next section) and we will find that a Bernoulli imitating process like the ones we build as $\Ber(\Pi)$ will contain enough of them if we take some care.

\section{Furcations}\label{s.furcations}

Furcation vertices will let us reach our target degree goal through Lemma \ref{l.furcationPlus2}, and will find trees containing them in Bernoulli clusters through Lemma \ref{l.basicBernoulli}.

\begin{definition}\label{d.Furcation} If $T$ is a tree, we say that  $v\in V(T)$ is a {\it furcation} if after deleting $v$ from $T$
among the remaining components there are at least $3$ infinite ones.
If a tree $T$ has a furcation we will say that $T$ is {\it forking}. When ${\cal F}$ is a forest and $v\in V({\cal F})$, then we will also say that $v$ is a furcation of ${\cal F}$ if it is a furcation of the subtree ${\cal F}(v)$ containing it. \end{definition}

As a trivial  example for a tree without furcation consider a bi-infinite path. The next example shows trees whose furcations are arranged in an adversarial  fashion and if $S$ is the set of their furcations, then the cells of ${\Vor}(S)$ are not finite. These examples are also worth keeping in mind, as they would be obstacles to our target degree goals. However, they simply cannot occur in an invariant process on a Cayley-graph (as we have seen in the application of the Mass Transport Principle).

\begin{example}\label{x.PathologicalFurcation} A {\it ray} emanating from $v$ is half-infinite path starting from $v$. Let $T_{\perp}$ be a tree (defined up to isomorphism) which  has a unique vertex of
 degree $3$ and all the other degrees are $2$ (i.e., three disjoint rays
 emanating from a single vertex).  As a further example consider the tree $T_{\perp\perp}$ obtained from a bi-infinite path $ P$ by attaching to each vertex $v$ a ray  $R_v$ (which are not interesecting each other or $P$). \end{example}

The following can be proven by induction on $r$.

\begin{lemma}\label{l.furcationPlus2} If $C$ is a finite, connected subset of a tree and contains at least $r$ furcations, then after deleting $C$ from the tree, among the remaining components there will be at least $ r+2$ which are infinite.\qed \end{lemma} 

\medskip
We mentioned that a $\Ber(p)$-labeling on $\T_d$ has infinite clusters whose label is constant $1$ if $p(d-1)>1$. This applies in particular to a $\Ber({1\over 2})$ labeling of $\T_4$. This implies that for a specific vertex $v$ of $\T_4$ the probability that $v$ will be in an infinite cluster in a $\Ber({1\over 2})$-labeling is positive. For us, however, forking clusters will be needed, and in the case of Bernoulli percolations,, we can easily get their existence from the infinite ones as follows.

For all $\epsilon>0$ there exists a $D(\epsilon)\in\N$,
 such that if $T$ is a tree whose minimal degree is at least $D(\epsilon)$ and $r$ is a distinguished vertex (the ``root'') of $T$, then in a  ${\Ber}({1\over 2})$-labeling of $T$ the cluster of the root $r$ will
 be a forking one with probability at least $1-\epsilon$. 
 
 To see this, consider first the rooted tree $(T,r)$ built from rooted copies $(T_1,r_1),\dots, (T_D,r_D)$ of $\T_4$ (so $T_i$ is a $4$-regular tree and $r_i$ is one of its vertices and $D$ will be fixed later) by adding a new vertex $r$ to this collection and make it into a tree by connecting $r$ to $r_i$ for all $i$ (no other edges are added). 
  For $r$ to be in a forking cluster in a $\Ber({1\over 2})$-labeling of $T$ it is enough if there are at least three such $(T_i,r_i)$ so that $r_i$ in an infinite cluster of the labeling restricted to $T_i$ and $r$ is connected to $r_i$. The probability of this clearly goes to $1$ as $D\rightarrow \infty$. If $D:=D(\epsilon)$ is chosen so that this probability is at least $1-\epsilon$, then in a tree whose minimal degree is at least $D$ we can take any vertex to be the root and embed this $(T,r)$ graph into it.

 Choose a sequence $\epsilon_0\in (0,1), \dots, \epsilon_n\in (0,1),\dots $
 tendig to $0$ fast enough so that $\prod_{n=1}^{\infty}
 (1-\epsilon_n)>0$. Define $D_n:=D(\epsilon_n)$ and use this sequence as the target degree in our construction.
 
\begin{lemma}\label{l.basicBernoulli} If a sequence $(T_1,r_1),\dots,(T_n,r_n),\dots $ of rooted
 trees is given, where the minimal degree of $T_n$ is at least $D_n$, and each of these trees are independently ${\Ber}({1\over 2})$-labeled, 
then with positive probability the roots of all
 of these trees will be in a forking cluster simultaneously. \end{lemma}
 
 We add a more process-oriented corollary to this. Assume that we start with a rooted tree $(T_0,r_0)$ which has minimal degree at least $D_0$, and we run the following process: label $T_0$ by a $\Ber({1 \over 2})$-labeling ${\cal L}_0$, and if the cluster of $r_0$ in $\clust({\cal L}_0)$ is not-forking, then stop, otherwise generate a new random rooted tree $(T_1,r_1)$ whose minimal degree is at least $D_1$. If the process has not stopped after $n$ steps, then we will have a rooted $(T_n,r_n)$ tree with minimal degree at least $D_n$. From here the process continues as in the beginning: let ${\cal L}_n$ be a $\Ber({1\over 2})$-labeling of $T_n$ (which is conditioned on $(T_n,r_n)$ is independent of the previous steps) and stop if the cluster of $r_n$ is not-forking, otherwise generate a random $(T_{n+1},r_{n+1})$ tree whose minimal degree is at least $D_{n+1}$. Then with positive probability the above process never stops and the generated sequence of bits ${\cal L}_0(r_0),\dots,{\cal L}_n(r_n),\dots$ will be i.i.d., so the random real number whose bits are this sequence has distribution $\un[0,1]$.
 
 \begin{remark} For us the distinction between an infinite and a forking tree is very important. However, it is known from \cite{Oded-Russ} that if $\lambda$ is {\it any} Bernoulli percolation on {\it any} Cayley graph, then its infinite clusters are {\it indistinguishable} by any invariant Borel property (which includes the one of being forking). Thus as soon as there are forking clusters with positive probability, then we know than in fact all infinite clusters are forking. Our labelings are not immediately Bernoulli ones, but cooked up from them in a way that this theorem would likely go through. However we did not try to use this direction as what we need can be obtained directly from the very basics of percolation theory on a tree.
 \end{remark}
 
 In our construction, Bernoulli processes on large scale forests (from Definition \ref{d.LargeScale}) will be used, and with the aid of the above, we will find forking forking ones among its clusters. We will put those clusters into use through the following lemma.

\begin{lemma}\label{l.DegGrowth} If $\Pi$ is an $\FIID$ cell-partition of $\T_3$ and $\F$ is an $\FIID$ subforest of $\T_3$ in such a way that $\Pi\prec \F$, and $D$ is a positive integer, then there is cell partition ${\tt {furc}}_D(\Pi)$ such that $\Pi\prec {\tt {furc}}_D(\Pi)\prec \F$ also holds, and whenever for a vertex $v$ the tree $\F(v)$ is forking, then the ${\tt {furc}}_D(\Pi)(v)$-cell contains at least $D$ furcations of $\F(v)$. \end{lemma}

\proof We first show that there exists a partition $\Pi^{\exists {\tt f}}$ for which $\Pi\prec\Pi^{\exists {\tt f}}\prec{\F}$ and whenever ${\F}(v)$ is forking, then $\Pi^{\exists {\tt f}}(v)$ contains at least one furcation of ${\F}(v)$.  Let ${\tt F}_{\Pi}$ be the set of those $\Pi$-cells which contain at least one furcation of ${\F}$. Our goal is achieved if we manage to glue $\Pi$-cells within a forking component to form bigger (but still finite) cells in such a way that every new cell contains at least one ``old'' $\Pi$-cell from ${\tt F}_{\Pi}$.

Voronoi cells on the large scale forest are just right for this purpose.
Move to the large scale forest ${\F}/\Pi$ and build $\Vor({\tt F}_{\Pi})$. This $\Vor({\tt F}_{\Pi})$ is ``almost'' the partition $\Pi^{\exists {\tt f}}$ we seek, except that it lives in ${\F}/\Pi$ instead of ${\F}$. We bring it back to ${\F}$ in the obvious way as $\Pi^{\exists {\tt f}}:={\tt glue}_{\Vor({\tt F}_{\Pi})}(\Pi)$.

Note that the finiteness of the new cells are guaranteed by induction, as every new cell either contains the distinguished finite subset which was an old cell from ${\tt F}_{\Pi}$, or (in case the ${F}$-component of a cell does not contain any furcation) the new cell is just equal to the old one.

Now that we have $\Pi^{\exists {\tt f}}$, we can define ${{\tt f}_D }(\Pi)$.   Since every $\Pi^{\exists {\tt f}}$-cell within a forking ${\F}$-cluster contains at least one furcation, it is enough if we manage to glue together $\Pi^{\exists {\tt f}}$-cells in such a way that every new cell of a forking cluster contains at least $D$ ``old'' $\Pi^{\exists {\tt f}}$-cells. To achieve this we can use the same idea as in the very first step described in the high level overview in constructing $\Pi_0$ and the associated large scale tree $\T_3/\Pi_0$, but this time we work within the forking components of ${\F}/\Pi^{\exists {\tt f}}$. In the large scale forest ${\F}/\Pi^{\exists {\tt f}}$ every forking component has minimal degree at least $3$. 
In ${\F}/\Pi^{\exists {\tt f}}$ select an $\FIID$ vertex set $S$ where the minimal distance between distinct vertices is at least $2D+1$ and $S$ has at least one element in every forking component of ${\F}/\Pi^{\exists {\tt f}}$.  In the corresponding Voronoi partition $\Vor(S)$,  every cell ${\cal C}$ belonging to a forking component of ${\F}/\Pi^{\exists {\tt f}}$ will contain at least $|B_{\T_3}(o,D)|\geq D$ many vertices of ${\F}/\Pi^{\exists {\tt f}}$ ($o$ denotes a generic vertex of $\T_3$). Thus we can define ${{\tt f}_D }(\Pi):={\tt glue}_{\Vor(S)}(\Pi^{\exists {\tt f}})$. The new cells are finite again by induction. \qed

\section{The main construction}\label{s.Main}

Now we give the precise definition of the sequence of cell-partitions $\Pi_0\prec\Pi_1\prec \dots\prec\Pi_n\prec \dots$; this will gives us also the sequence of $\Lambda_i\sim\Ber(\Pi_i)$ labels, where conditioned on $\Pi_i$ the label $\Lambda_i$ will be independent from the previous labels (but $\Pi_i$ itself depends on  $\{(\Pi_j,\Lambda_j)\}_{j<i}$).

$\Pi_0$ and $\Lambda_0$ are as defined before in~Section \ref{s.High_Level}. Assume that $\Pi_0\prec \dots \prec\Pi_n$ and $\Lambda_0,\dots,\Lambda_n$ are defined. Let ${\cal F}_n:=\clust(\Lambda_0,\dots,\Lambda_n)$, where $(\Lambda_0,\dots,\Lambda_n)$ is the $\{0,1\}^{n+1}$-label obtained by concatenating the $\Lambda_i$s.

We want to define $\Pi_{n+1}$ in such a way that $\Pi_n\prec \Pi_{n+1} \prec {\cal F}_n$, and if for a vertex $v$ the tree ${\cal F}_n(v)$ is forking, then the $\Pi_{n+1}(v)$-cell should contain at least $D_{n+1}$ furcations of ${\cal F}_n$. 
We use Lemma \ref{l.DegGrowth} for the pair $\Pi_n\prec {\cal F}_n$ and define $\Pi_{n+1}:={\tt f}_{D_{n+1}}(\Pi_n)$.

This concludes the construction of $\Pi_0\prec\dots\prec\Pi_n\prec\dots$ and thus also that of $\Lambda_i\sim\Ber(\Pi_i)$, with the specification that conditioned on $\Pi_i$ the $\Lambda_i$ must be independent of the previous steps (which implies that for a generic vertex $o$ the sequence $\Lambda_0(o),\Lambda_1(o),\dots,$ of bits is i.i.d.).
Because $\Pi_{n+1}\prec{\cal F}_n:=\clust(\Lambda_0,\dots,\Lambda_n)$, if we define $\Pi_{\infty}(o):=\bigcup_{i\in \N}\Pi_i(o)$, then for any $v_1,v_2\in \Pi_{\infty}(o)$ and $m\in \N$ we have $\Lambda_m(v_1)=\Lambda_m(v_2)$.
Thus if we define $\Lambda_{\infty}(o)$ to be the real number from $[0,1]$ whose consecutive bits are $\Lambda_0,\dots,\Lambda_n,\dots$ (which is an i.i.d.\ sequence, so $\Lambda_{\infty}(o)$ has distribution $\un[0,1]$), then $v_1,v_2\in \Pi_{\infty}(o)$ also implies $\Lambda_{\infty}(v_1)=\Lambda_{\infty}(v_2)$, thus $\Pi_{\infty}(o)$ will be contained within a single cluster of $\clust(\Lambda_{\infty})$.
Moreover, if ${\cal F}_n(o)$ is a forking cluster, then $|\Pi_{n+1}(o)|\geq D_{n+1}$ as it contains at least $D_{n+1}$ furcation of ${\cal F}_n(o)$. So if ${\cal F}_i(o)$ is forking for every $i$, then $\Pi_{\infty}(o)\supset \Pi_m(o)$ contains at least $D_m$ element for any $m$, and as $D_m\to \infty$, this implies ~Theorem\ref{t.gab}.

To conclude, we observe that indeed it happens with positive probability that ${\cal F}_n(o)$ is forking for all $n$, because of the corollary to ~Lemma \ref{l.basicBernoulli} using the process of generating rooted trees. The correspondence is as follows. Start with the rooted tree $(T_0,r_0):=(\T_3/\Pi_0(o),\Pi(o))$ whose minimal degree is greater than $D_0$. Use our $\Lambda_0\sim\Ber(\Pi_0)$ which is a $\Ber({1\over 2})$-labeling of $T_0$. Stop the process if the cluster ${\cal F}_0(o)$ of $\Pi_0(o)$ is not forking, otherwise continue by creating the next random rooted tree $({\cal F}_0/\Pi_1(o),\Pi_1(o))$, whose minimal degree is greater than $D_1$. In general, if the process has not stopped, then the rooted tree  $(T_n, r_n)$ is constructed as the random  rooted tree $({\cal F}_n/\Pi_{n+1}(o),\Pi_{n+1}(o))$. Notice that -- by the finiteness of the cells -- moving from the tree ${\cal F}_n(o)$ to its large scale version ${\cal F}_n/\Pi_{n+1}(o)$ does not change its being forking or not (while it increases its minimal degree), so the correspondence between our construction and the process oriented corollary to ~Lemma \ref{l.basicBernoulli} is complete. Our example is also manifestly hyperfinite because not only $\Pi_{\infty}(o)\subset \clust(\Lambda_{\infty})$ for any $o$ which already means $\Pi_{\infty}\prec\clust(\Lambda_{\infty})$ but actually $\Pi_{\infty}=\clust(\Lambda_{\infty})$ (and $\Pi_{\infty}(o)$ is an increasing union of the finite $\Pi_n(o)$s). To see this, consider an edge $e$ connecting a vertex $v_{\tt in}\in \Pi_{\infty}(o)$ with $v_{\tt out}\not\in \Pi_{\infty}(o)$, notice that for any $i$, the labels $\Lambda_i(v_{\tt in})$ and $\Lambda_i(v_{\tt out})$ are independent, so they cannot be all equal. So $e$ is deleted from $\clust(\Lambda_{\infty})$.

\section{Merging small clusters into big ones and the second construction}\label{s.merging}

Assume that we classify trees as ``small" and ``big" and our
classification scheme has the natural property that being big is
upward closed in the sense that if $T$ is a tree which contains a big
subtree, then $T$ is itself big. Two natural examples are: ``being
infinite" and ``being forking".
We show that if we have a $\{0,1\}$-labeling which
almost surely has big clusters, then there is a natural way to
``merge" the small clusters into the big ones so that at the end only
big ones remain. It will be achieved through a
``relabeling" which replaces the old label ${\cal L}$ by a new one
${\cal L}^*$ so that all the clusters of the ${\cal L}^*$ label are
big.

This merging process will not need any extra randomness. It consist of iterating the following:
 if a small ${\cal
L}$-cluster ${ C}$ is at distance $1$ from a big one, then
every vertex in ${ C}$ switches its label so that ${ C}$
``joins" the forking cluster. Note that there is no ambiguity: if there are more than one big clusters at distance $1$, then each must have the same label as we only have two labels.
In this way the big clusters have grown and we repeat the process. The fact that any small cluster is at finite distance from some big one (simply by the existence of big clusters) implies that every cluster will join a big one eventually.

Using this relabeling we can modify the construction so that every cluster will be infinite.

 This is also built through a sequence of cell partition $P_0,\dots,P_n,\dots$, but the labeling will not be simply a version of $\Ber(P_i)$. We assume that $P_0,\dots,P_n$ and labels $L_0,\dots,L_n$ are already defined. As before, this gives $F_n:=\clust((L_0,\dots,L_n))$, and it will also be true by the iteration step that $P_n\prec F_n$. Moreover, we will have that in the large scale forest $F_n/P_n$ {\it every} component is forking.
 We define $P_{n+1}$ in such a such a way that in the large scale forest $F_n/P_{n+1}${\it every} component has minimal degree at least $4$ (this can be done just as in the previous section with the target degree being $4$ this time). Then first define the labeling $l\sim\Ber(P_{n+1})$ (independent of the previous steps conditioned on $P_{n+1}$), then $\clust(l)$ will contain forking clusters so we can use the relabeling construction to get $l^*=:L_{n+1}$. This gives a labeling as claimed.

\section{Acknowledgements.} This work was partially supported by the ERC Consolidator Grant 772466 ``NOISE''. At earlier stages the work was also partially supported by NSF Grant DMS-1007244, and OTKA Grant K76099. I am grateful for Russell Lyons for suggesting working on this problem and for useful discussions and for helping improving the first version of the text. I also thank the anonymous referee of the first submitted version for useful feedback.
I am indebted to Gábor Pete for useful discussions throughout the writing of this version and for Sándor Rokob for useful feedback.

\end{document}